\documentclass[graybox]{svmult}

\usepackage{mathptmx}       
\usepackage{helvet}         
\usepackage{courier}        
\usepackage{type1cm}        
\usepackage{makeidx}         
\usepackage{graphicx}        
\usepackage{multicol}        
\usepackage[bottom]{footmisc}

\usepackage{cite}
\usepackage{fancyvrb}        

\usepackage{amsmath}
\usepackage{amssymb}
\usepackage{bm}

\makeindex             

\begin{document}

\title*{Wavelet Algorithm for Circuit Simulation}
\author{Kai Bittner \and Emira Dautbegovic }
\institute{Kai Bittner\at University of Wuppertal, Germany,
  \email{bittner@math.uni-wuppertal.de}
  \and Emira Dautbegovic \at Infineon Technologies AG, 81726 Munich, Germany,
  \email{Emira.Dautbegovic@infineon.com}
}
%
%
\maketitle

\abstract*{We present a new adaptive circuit simulation algorithm based on 
spline
wavelets. The unknown voltages and currents are expanded into a wavelet 
representation, which is determined as solution of nonlinear equations derived
from the circuit equations by a Galerkin discretization. The spline wavelet
representation is adaptively refined during the Newton iteration.
The resulting approximation requires an almost minimal number of degrees of 
freedom, and in addition the grid refinement approach enables very efficient 
numerical computations. Initial numerical tests on various
types of electronic circuits show promising results when compared to the 
standard transient analysis.} 

\abstract{We present a new adaptive circuit simulation algorithm based on 
spline
wavelets. The unknown voltages and currents are expanded into a wavelet 
representation, which is determined as solution of nonlinear equations derived
from the circuit equations by a Galerkin discretization. The spline wavelet
representation is adaptively refined during the Newton iteration.
The resulting approximation requires an almost minimal number of degrees of 
freedom, and in addition the grid refinement approach enables very efficient 
numerical computations. Initial numerical tests on various
types of electronic circuits show promising results when compared to the 
standard transient analysis.} 

\section{Introduction}

Wavelet theory emerged during the $20^\mathrm{th}$ century from the
study of Calderon-Zygmund operators in mathematics, the study of the
theory of subband coding in engineering and the study of
renormalization group theory in physics. Recent approaches
\cite{p208:BaKnPu08,p208:CS01,p208:DCB05,p208:SGN07,p208:ZC99} 
to the problem of multirate
envelope simulation indicate that wavelets could also be used to
address the qualitative simulation challenge by a development of novel
wavelet-based circuit simulation techniques capable of an efficient
simulation of a mixed analog-digital circuit \cite{p208:ED_SCEE08}.

The wavelet expansion of a function $f$ is given as
\begin{equation}\label{p208:wavelet_expansion}
f= \sum_{k\in\mathcal{I}} c_k\, \phi_k + 
\sum_{j=0~}^\infty \sum_{k\in\mathcal{K}_j}d_{jk}\,\psi_{jk}.
\end{equation}
Here, $j$ refers to a level of resolution, while $k$ describes the
localization in time or space, i.e., $\psi_{jk}$ is essentially
supported in the neighborhood of a point $x_{jk}$ of the domain,
 where the wavelet expansion
is defined. The wavelet
expansion can be seen as coarse scale approximation 
$\sum_{k\in\mathcal{I}} c_k\, \phi_k$ by the scaling functions $\phi_k$
complemented by information on details of increasing resolution $j$
in terms of the wavelets $\psi_{jk}$.
Since a wavelet basis consist of an infinite number of wavelets one
has to consider approximations of $f$ by partial sums of the wavelet
expansion (\ref{p208:wavelet_expansion}), which can, e.g.,  be obtained
by ignoring small coefficients.  

\section{Wavelet-based Galerkin Method}

We consider circuit equations in the charge/flux oriented 
modified nodal analysis (MNA) formulation, which yields a mathematical model 
in the form of an initial-value problem of differential-algebraic equations 
(DAEs):  
\begin{equation}
  \label{p208:eq_MNA_charge}
  \frac{d}{dt}\bm{q}\big(\bm{x}(t)\big) 
        + \bm{f}\big(\bm{x}(t)\big) = \bm{s}(t). 
\end{equation}  
Here 
$\bm{x}$ is the vector of node potentials and specific branch voltages and 
$\bm{q}$ is the vector of charges and fluxes. Vector $\bm{f}$ comprises static 
contributions, while $\bm{s}$ contains the contributions of independent
sources. 

In our adaptive wavelet approach we first discretize the MNA equation 
(\ref{p208:eq_MNA_charge}) in terms of the wavelet basis functions, by 
expanding $\bm{x}$ as a linear combination of wavelets
or related basis functions $\varphi_i$, i.e., $\bm{x}=\sum_{i=0}^n \bm{c}_i\,\varphi_i$. 
Furthermore,  we integrate  the circuit equations against test functions
 $\theta_\ell$ and obtain the equations 
\begin{equation}
\label{p208:nonlinear}
\int_0^T \Big(\frac{d}{dt}\bm{q}\big(\bm{x}(t)\big) 
        + \bm{f}\big(\bm{x}(t)\big) - \bm{s}(t)\Big)\,\theta_\ell\,dt = 0,
\end{equation} 
for $\ell=1,\ldots,n$.
Together with the initial conditions $\bm{x}(0)= \bm{x}_0$, we have now
$n+1$ vector valued equations, which determine the coefficients $\bm{c}_i$
provided that the test functions
$\theta_\ell$ are chosen suitably to the basis functions
$\varphi_i$. 

The nonlinear system (\ref{p208:nonlinear}) is solved by Newton's method.
With a good initial guess, Newton's method is known to converge quadratically.
Unfortunately, a good initial guess is usually not available, and convergence 
can often only be obtained after a large number of (possibly damped) Newton steps. 
On the other
hand, to get a good approximation of the solution of 
(\ref{p208:eq_MNA_charge}), the space $X =\mathrm{span}\{\varphi_k:~k=0,\ldots,n\}$
has to be sufficiently large and the computational cost of each step depends
on $n=\dim X$.

Here, we take advantage from the use of wavelets. The Newton iteration is started
on a coarse subspace $X_0$ of small dimension, which provides us with a first coarse
approximation $x^{(0)}$ of the solution. Then $x^{(0)}$ is used as initial
guess for 
a Newton iteration in a finer space $X_1 \supset X_0$, leading to an improved approximation
$x^{(1)}$. One positive effect, which can be observed in numerical tests, is that
a single Newton step in the beginning of the algorithm is relatively cheap, i.e.,
having only a poor initial guess for $x^{(i)}$ with $i$ small has only a negligible effect
on the performance of the algorithm. On the other hand, due to the excellent initial guess
in the higher dimensional spaces $X_i$ with $i$ large, we need only a few of the costly Newton
steps, which are necessary in order to achieve a required accuracy.
The embedding $X_i \subset X_{i+1}$ is ensured by the use of wavelet subsets, i.e.,
$$
X_i = \mathrm{span}\Big(\{\phi_k:k\in\mathcal{I}\}\cup \{\psi_{jk}:~(j,k)\in\Lambda_i\}\Big),\qquad
\Lambda_i\subset \Lambda_{i+1},
$$
i.e., we add adaptively more and more wavelets to the expansion.


Due to the intrinsic properties of wavelets \cite{p208:ED_SCEE08} an adaptive wavelet 
approximation can provide an efficient representation
of functions with steep transients, which often appear in a mixed 
analog/digital electronic circuit.
However, for an efficient circuit simulation we have to take
further properties of a wavelet system into account. We consider spline
wavelets to be the optimal choice since 
spline wavelets are the only wavelets with an explicit formulation.
This permits the fast computation of function values, derivatives and integrals,
which is essential for an efficient solution of a nonlinear problem as given in 
(\ref{p208:eq_MNA_charge}) (see also \cite{p208:BiUr04, p208:DSX00}).
Spline wavelets have already been used for circuit simulation \cite{p208:ZhCa99}.
However, here we use a completely new approach based on spline wavelets from
\cite{p208:Bit05b}.

\section{Numerical Tests \label{p208:numtest}}

A prototype of the proposed wavelet algorithm 
is implemented within the framework of a productively used circuit 
simulator and tested on a variety of circuits.
For all examples we have compared the CPU time
and the grid size (i.e., the number of spline knots or time steps) with
the corresponding results from transient analysis of the same
circuit simulator. 

The error is estimated by comparison with well-established and
highly-accurate transient analysis. The estimate shown in the signal is the 
maximal absolute difference over all grid points of the transient analysis, which gives a good
approximation of the maximal error. That is, if we can obtain a small
error for the wavelet analysis, which proves good agreement with the standard
method. In particular, since we compare the solutions of two independent methods
we have very good evidence that we approximate the solution of the underlying 
DAE's with the estimated error.

\subsection{Amplifier}

\begin{figure}[htbp]
\centering
\begin{tabular}{c}
\includegraphics[width=0.7\textwidth]{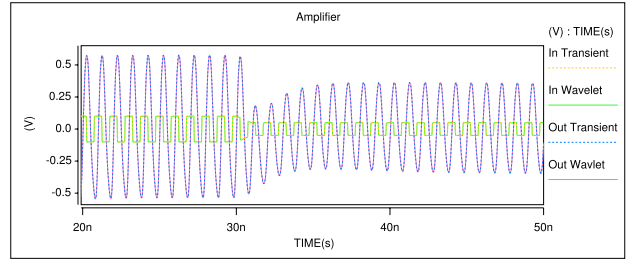}
\end{tabular}
\caption{Detail of Input and output signal for the amplifier \label{p208:amp_out_detail}}
\end{figure}

The amplifier is simulated with a pulse signal of period 1 ns, which
is modulated by a piecewise smooth amplitude 
(see Fig.~\ref{p208:amp_out_detail}). 
The wavelet method runs over 100 ns. 
The results show a satisfying performance also for digital-like
input signals.

\begin{figure}[htbp]
\centering
\begin{tabular}{cc}
\includegraphics[width=0.45\textwidth]{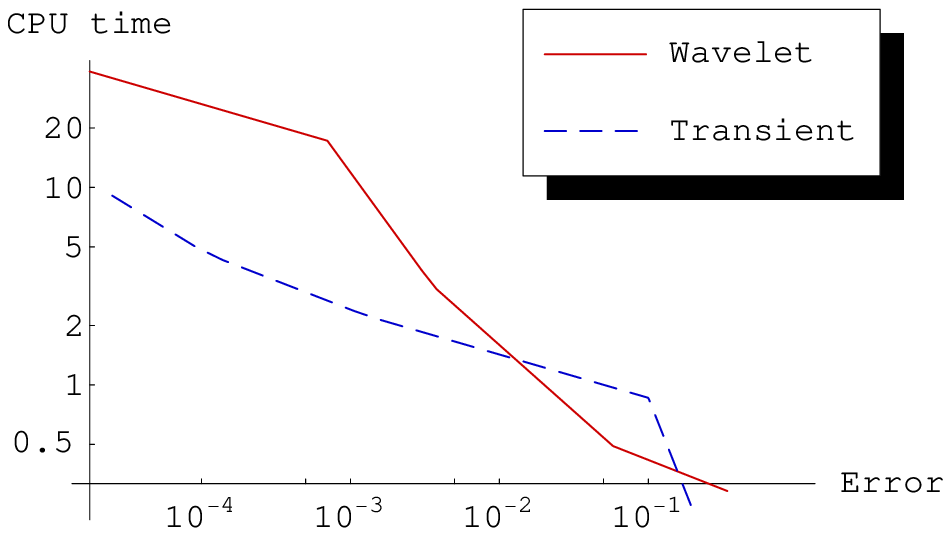}
 &\includegraphics[width=0.45\textwidth]{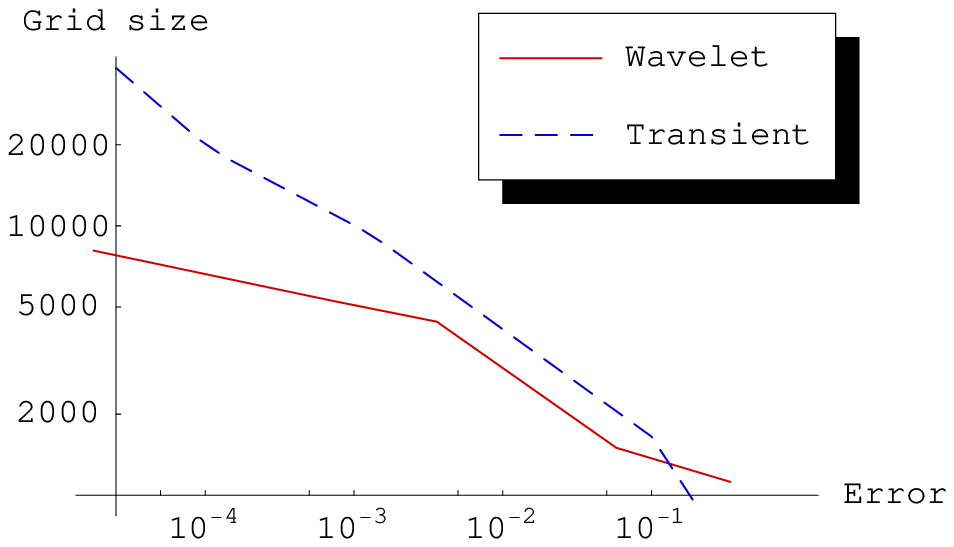}
\end{tabular}
\caption{Simulation statistics for the amplifier.
Computation time versus error (left), and
grid size versus error (right) for transient analysis and adaptive
wavelet analysis.
\label{p208:amp_res}}
\end{figure}

\subsection{Oscillator}

\begin{figure}[htbp]
\centering
\begin{tabular}{cc}
\includegraphics[width=0.45\textwidth]{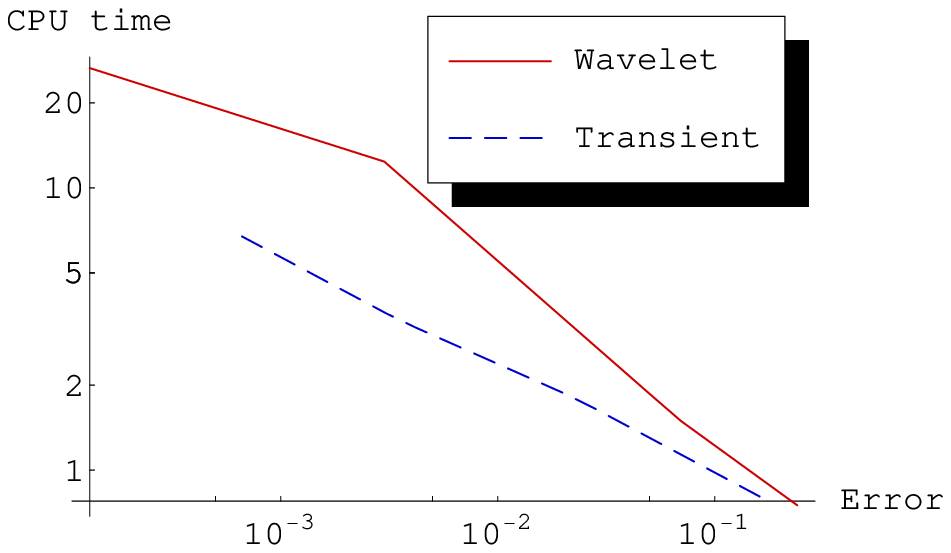}
 &\includegraphics[width=0.45\textwidth]{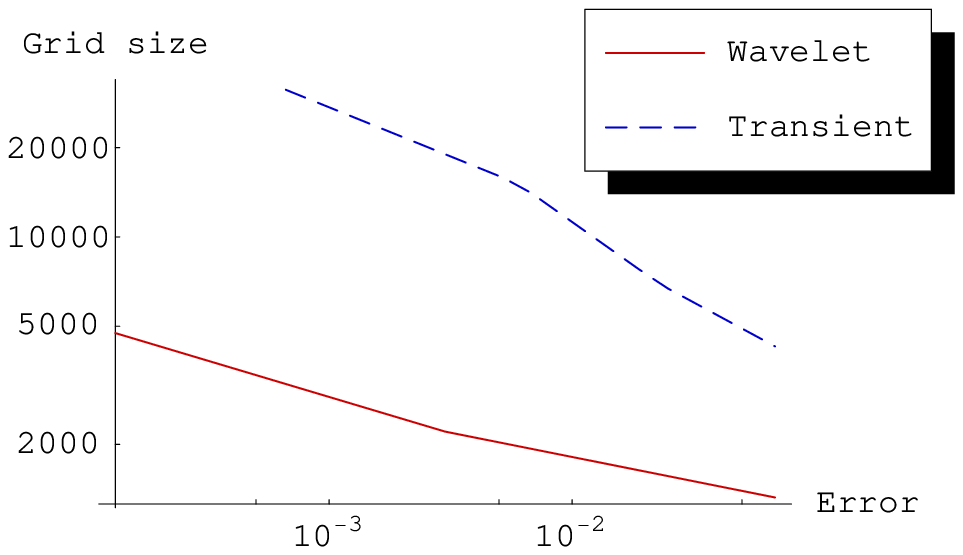}
\end{tabular}
\caption{Simulation statistics for the oscillator.
Computation time versus error (left), and
grid size versus error (right) for transient analysis and adaptive
wavelet analysis.
\label{p208:vcobi_res}}
\end{figure}

The oscillator is an autonomous circuit without an external input signal. 
The simulation runs over 20 ns. As can be seen from Fig~\ref{p208:vcobi_res},
an excellent agreement with highly-accurate transient analysis is achieved. 
It should be noted that after the oscillator has reached its periodic 
steady state the wavelet method works very fast, since the solution from 
one interval is an excellent initial guess for the next interval. 

\begin{figure}[htbp]
\centering
\begin{tabular}{c}
\includegraphics[width=0.7\textwidth]{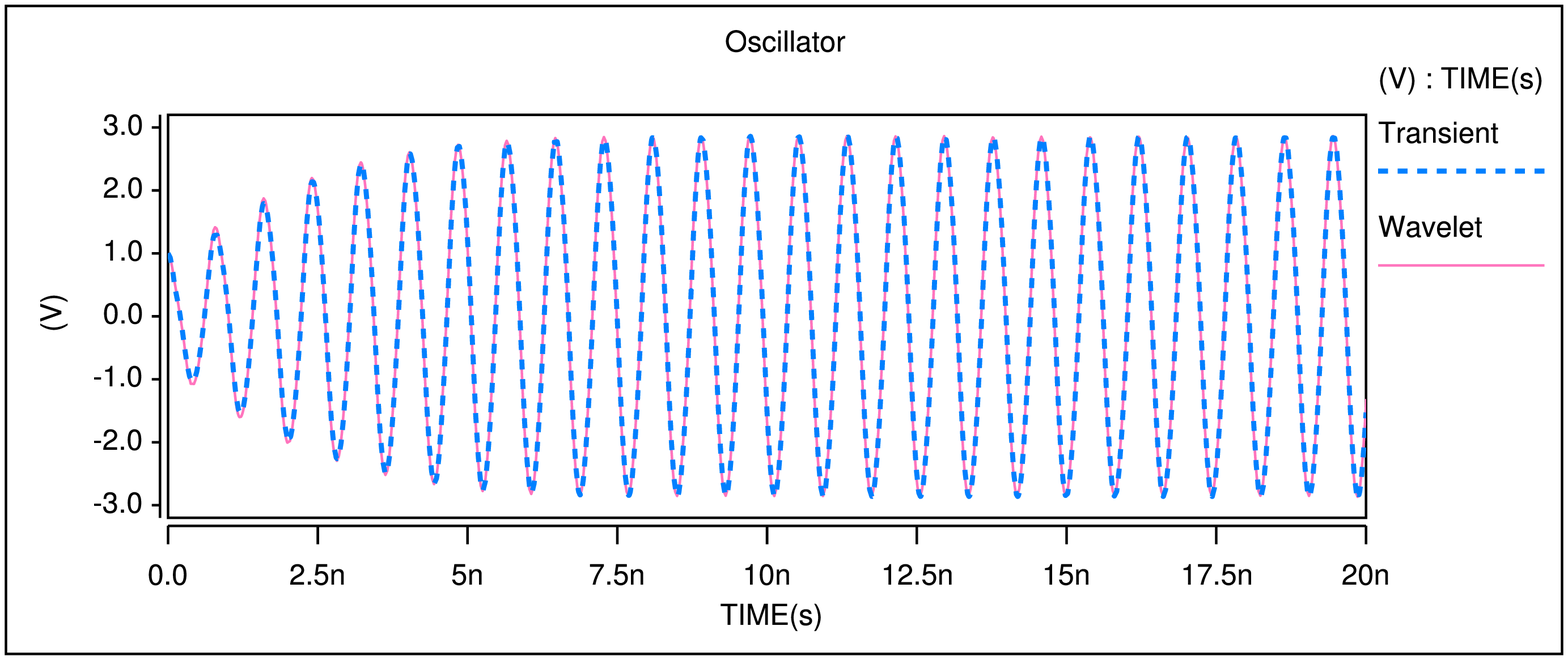}
\end{tabular}
\caption{Output signal for the oscillator.
\label{p208:vcobi_out}}
\end{figure}

\subsection{Mixer}

The mixer is simulated with input frequencies 
950 MHz and 1GHz.
The simulation runs over 30 ns. 
In particular, for high 
accuracies the number of degrees of freedom is essentially reduced, while
the computation time is at least of the same order. 

\begin{figure}[htbp]
\centering
\begin{tabular}{c}
\includegraphics[width=0.7\textwidth]{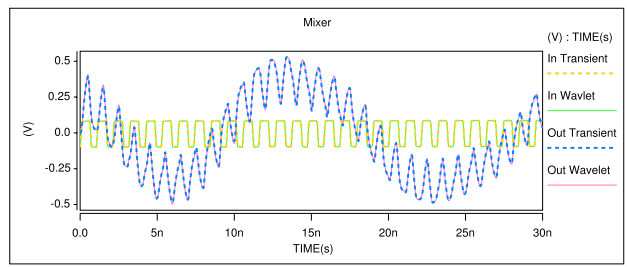}
\end{tabular}
\caption{Time domain output signal for the mixer.
\label{p208:mixer_out}}
\end{figure}

\begin{figure}[htbp]
\centering
\begin{tabular}{cc}
\includegraphics[width=0.45\textwidth]{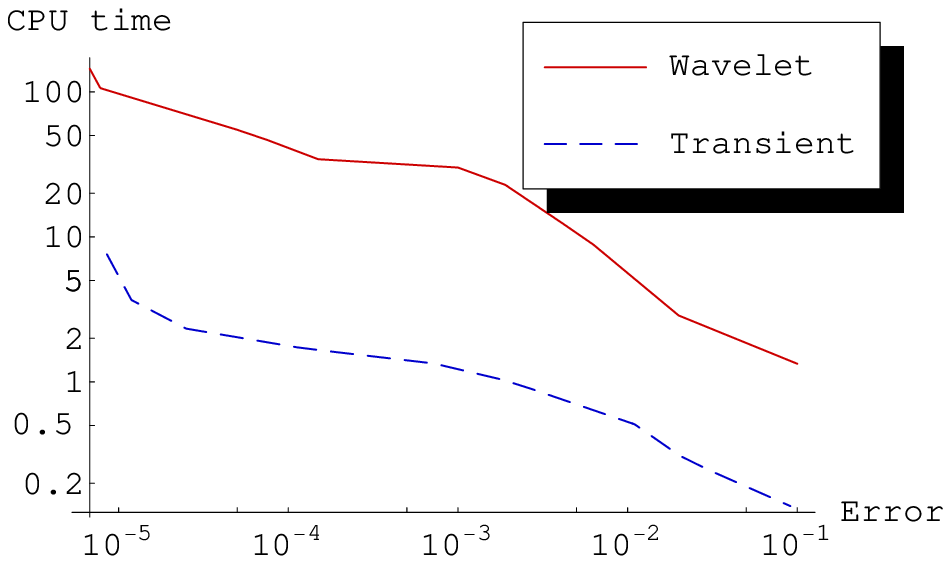}
 &\includegraphics[width=0.45\textwidth]{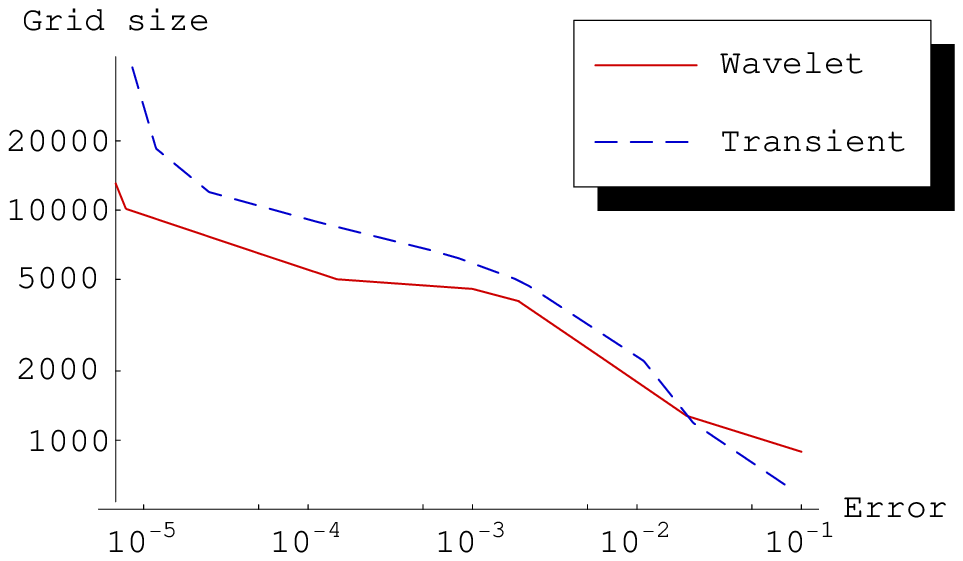}
\end{tabular}
\caption{Simulation statistics for the mixer.
Computation time versus error (left), and
grid size versus error (right) for 
transient analysis and adaptive
wavelet analysis.
\label{p208:mixer_res}}
\end{figure}

\section{Conclusion}

The results of the simulations indicate that the wavelet based method
may achieve and in some cases surpass performance of the standard transient 
analysis. 
Apparently, the number of degrees of freedom can be smaller than for the 
transient analysis for comparable accuracy. However, this advantage of 
the wavelet algorithm does not always result (yet) in a smaller
computation time. 
On the other hand it can be expected that the productive implementation 
of the wavelet algorithm can be further optimized. Therefore our activities on 
optimization and further development of the wavelet-based algorithm are continuing.  

\begin{acknowledgement}
  This work has been supported within the EU Seventh Research Framework 
  Project (FP7) \mbox{ICESTARS} with the grant number 214911.  
\end{acknowledgement}

%
%
%




\end{document}